\documentclass[11pt]{amsart}
\usepackage{amsmath,amsthm,amsfonts,amssymb,latexsym,epsfig}

\newtheorem{theorem}{Theorem}[section]

\newtheorem{lemma}{Lemma}[section]

\numberwithin{equation}{section}

\theoremstyle{definition}

\theoremstyle{remark}

\begin{document}
\title{On $l^2$ norms of some weighted mean matrices}
\author{Peng Gao}
\address{Department of Computer and Mathematical Sciences,
University of Toronto at Scarborough, 1265 Military Trail, Toronto
Ontario, Canada M1C 1A4} \email{penggao@utsc.utoronto.ca}
\date{March 10, 2008.}


\begin{abstract}
 We give another proof of a result of Bennett on the $l^{p}$ operator norms of some weighted mean matrices for the case $p=2$ and we also present some related results.
\end{abstract}

\maketitle
\section{Introduction}
\label{sec 1} \setcounter{equation}{0}

  Suppose throughout that $p >1, \frac{1}{p}+\frac{1}{q}=1$.
   Let $l^p$ be the Banach space of all complex sequences ${\bf a}=(a_n)_{n \geq 1}$ with norm
\begin{equation*}
   ||{\bf a}||_p: =(\sum_{n=1}^{\infty}|a_n|^p)^{1/p} < \infty.
\end{equation*}
  The celebrated
   Hardy's inequality (\cite[Theorem 326]{HLP}) asserts that for $p>1$,
\begin{equation}
\label{eq:1} \sum^{\infty}_{n=1}\Big{|}\frac {1}{n}
\sum^n_{k=1}a_k\Big{|}^p \leq \Big (\frac
{p}{p-1} \Big )^p\sum^\infty_{k=1}|a_k|^p.
\end{equation}
   Hardy's inequality can be regarded as a special case of the
   following inequality:
\begin{equation*}
\label{01}
   \sum^{\infty}_{j=1}\big{|}\sum^{\infty}_{k=1}c_{j,k}a_k
   \big{|}^p \leq U \sum^{\infty}_{k=1}|a_k|^p,
\end{equation*}
   in which $C=(c_{j,k})$ and the parameter $p$ are assumed
   fixed ($p>1$), and the estimate is to hold for all complex
   sequences ${\bf a}$. The $l^{p}$ operator norm of $C$ is
   then defined as the $p$-th root of the smallest value of the
   constant $U$:
\begin{equation*}
\label{02}
    ||C||_{p,p}=U^{\frac {1}{p}}.
\end{equation*}

    Hardy's inequality thus asserts that the Ces\'aro matrix
    operator $C$, given by $c_{j,k}=1/j , k\leq j$ and $0$
    otherwise, is bounded on {\it $l^p$} and has norm $\leq
    p/(p-1)$. (The norm is in fact $p/(p-1)$.)

    We say a matrix $A$ is a summability matrix if its entries satisfy:
    $a_{j,k} \geq 0$, $a_{j,k}=0$ for $k>j$ and
    $\sum^j_{k=1}a_{j,k}=1$. We say a summability matrix $A$ is a weighted
    mean matrix if its entries satisfy:
\begin{equation*}
    a_{j,k}=\lambda_k/\Lambda_j,  ~~ 1 \leq k \leq
    j; \Lambda_j=\sum^j_{i=1}\lambda_i, \lambda_i \geq 0, \lambda_1>0.
\end{equation*}

    A natural generalization of Hardy's inequality \eqref{eq:1} is to
    determine the $l^{p}$ operator norm of an arbitrary summability matrix $A$.
 For examples, the following two
     inequalities were claimed to hold by Bennett ( \cite[p. 40-41]{B4}; see also \cite[p. 407]{B5}):
\begin{eqnarray}
\label{7}
   \sum^{\infty}_{n=1}\Big{|}\frac
1{n^{\alpha}}\sum^n_{i=1}(i^{\alpha}-(i-1)^{\alpha})a_i\Big{|}^p &
\leq & \Big( \frac {\alpha p}{\alpha p-1} \Big )^p\sum^{\infty}_{n=1}|a_n|^p, \\
\label{8}
   \sum^{\infty}_{n=1}\Big{|}\frac
1{\sum^n_{i=1}i^{\alpha-1}}\sum^n_{i=1}i^{\alpha-1}a_i\Big{|}^p &
\leq & \Big(\frac {\alpha p}{\alpha p-1} \Big
)^p\sum^{\infty}_{n=1}|a_n|^p,
\end{eqnarray}
     whenever $\alpha>0, p>1, \alpha p >1$. We note here the constant $(\alpha p /(\alpha p-1))^p$ is best possible (see \cite{Be1}).

     No proofs of the above two inequalities were supplied in \cite{B4}-\cite{B5}. The author \cite{G} and Bennett himself \cite{Be1}
     proved inequalities \eqref{7} for $p>1, \alpha \geq 1, \alpha p >1$ and
     \eqref{8} for $p>1, \alpha \geq 2$ or $0< \alpha \leq 1, \alpha p >1$
     independently. Bennett in fact was able
     to prove \eqref{7} for $p > 1, \alpha >0, \alpha p >1$ (see
     \cite[Theorem 1]{Be1} with $\beta=1$ there) which now leaves the
     case $p>1, 1< \alpha <2$ of inequality \eqref{8} the only case
     open. For this, Bennett expects inequality \eqref{8} to hold for $1+1/p < \alpha <2$
     (see page 830 of \cite{Be1}) and as a support,
     Bennett \cite[Theorem 18]{Be1} has shown that
inequality \eqref{8}
   holds for $\alpha = 1+1/p, p > 1$. Recently, much progress was made on this later case this by the author in \cite{G3} and \cite{G5}.

   It is our goal in this paper to study inequalities \eqref{7}-\eqref{8} for the case $p=2$, using the approach of quadratic forms.  For the case of Hardy's inequality \eqref{eq:1}, such an approach was used by Schur in \cite {Schur}, who showed that for ${\bf x}, {\bf y} \in l^2$, 
\begin{equation}
\label{1.5}
   \sum^{\infty}_{i,j=1}\frac {x_iy_j}{\max (i,j)} \leq 4 || {\bf x}||_2 || {\bf y}||_2. 
\end{equation}
   As we shall see later in this paper that the above inequality implies Hardy's inequality \eqref{eq:1}, even though this was not mentioned in \cite{Schur} (this is actually prior to Hardy's discovery of \eqref{eq:1}). Further developments were carried out by Wilf in \cite{W1} and Wang and Yuan in \cite{W&Y}. The idea of this approach via quadratic forms is to
   interpret the left-hand side of \eqref{eq:1} when $p=2$ as a quadratic form
   so that Hardy's inequality follows from an estimations of the largest
   eigenvalue of the corresponding matrix associated to the quadratic form.
   We will adopt this approach in Section \ref{sec 3} to study the $l^2$ norms of weighted mean matrices associated to \eqref{7} to give another proof of the following result:
\begin{theorem}
\label{thm1}
  Let $p=2$, then inequality \eqref{7} holds for $1/2 < \alpha \leq 3/2$ and inequality \eqref{8} holds for $1/2 < \alpha \leq 1$.
\end{theorem}

   In relation to Hardy's inequality \eqref{eq:1} for $p=2$, we note the following well-known Hilbert's inequality \cite[Theorem 315]{HLP} for ${\bf x}, {\bf y} \in l^2$:
\begin{equation}
\label{1.6}
   \sum^{\infty}_{i,j=1}\frac {x_iy_j}{i+j} \leq \pi || {\bf x}||_2 || {\bf y}||_2. 
\end{equation}
   One may regard, other than a constant factor, the entries of both the coefficient matrices of the quadratic forms of the left-hand side expressions of \eqref{1.5} and \eqref{1.6} as special cases of the following family of matrices for $r \geq 1$:
\begin{equation}
\label{1.7}
   \Big ( P^{-1}_r(i,j) \Big )_{i,j}, \hspace{0.1in} P_r(i,j)=\Big(\frac {i^r+j^r}{2} \Big )^{1/r}. 
\end{equation}
  Note here that the case $r=1$ above gives the case for \eqref{1.6}, except for a factor of $2$ and the case $r \rightarrow +\infty$ gives the corresponding case for \eqref{1.5}. We point out here the norms of the matrices given in \eqref{1.6}, or more generally, the norms of matrices whose $(i,j)$-th entry given by $K(i,j)$ with $K(x,y)$ positive, decreasing with respect to both $x$ and $y$ and homogeneous of degree $-1$ was studied by Schur in \cite{Schur}, see also \cite[Theorem 318-320]{HLP}.

  We note also that there is a different version of Hilbert's inequality \cite[Theorem 294]{HLP}, which asserts that for a sequence of complex numbers ${\bf w}$,
\begin{equation}
\label{1.8} 
   \Big |
\sum_{r \neq s} \frac {w_r\overline{w_s}}{r-s}\Big | \leq \pi \sum_{r}|w_r|^2.
\end{equation}
  The best constant $\pi$ was first determined by Schur \cite{Schur}. One may regard the entries of the coefficient matrix on the left-hand side expression of the above inequality as obtained from those of \eqref{1.6} by replacing the ``$+$'' sign by the ``$-$'' sign and omitting the diagonal entries. In our proof of Theorem \ref{thm1}, we shall make use a quadratic form whose coefficient matrix is essentially (other than a constant factor) given by the case of $r \rightarrow +\infty$ of the following family of matrices of $r \geq 1$:
\begin{equation}
\label{1.10}
   \Big ( P^{-1}_r(i^{\alpha}j^{1-\alpha},i^{1-\alpha}j^{\alpha}) \Big )_{i,j}, 
\end{equation}
  where $\alpha$ is a parameter and $P_r$ is given as in \eqref{1.7}. In view of our discussions above, it is natural to ask whether one can establish inequalities similar to \eqref{1.8} with the coefficient matrix given by the case $r=1$ of \eqref{1.10} with the ``$+$'' sign replaced by the ``$-$'' sign and omitting the diagonal entries. This is indeed possible and was done by Schur in \cite{Schur}. In what follows, we shall give a more general statement of Schur's result by first pointing out that Montgomery and Vaughan \cite{M&V} (see also \cite[Theorem 2]{M}) has given the following generalization of Hilbert's inequality \eqref{1.8}:
\begin{theorem}
\label{thm2'}
  Suppose that $\lambda_1<\lambda_2<\ldots<\lambda_R$, and that $\lambda_{r+1}-\lambda_r \geq \delta>0$ for $1 \leq r <R$. Then for any ${\bf w}$,
\begin{equation}
\label{2} 
   \Big |
\sum_{r \neq s} \frac {w_r\overline{w_s}}{\lambda_r-\lambda_s}\Big | \leq \pi \delta^{-1}\sum_{r}|w_r|^2.
\end{equation}
\end{theorem}
  The above inequality is useful in deriving the formulation of large sieves, which have applications in number theory. Schur's result can now be formulated in the following way:
\begin{theorem}
\label{thm2} 
  Let $\alpha \geq 1$ be fixed. Suppose that $\lambda_1<\lambda_2<\ldots<\lambda_R$, and that $\lambda_{r+1}-\lambda_r \geq \delta>0$ for $1 \leq r <R$. Then for any ${\bf w}$,
\begin{equation*}
   \Big |
\sum_{r \neq s} \frac {(\lambda_r\lambda_s)^{\frac {\alpha-1}{2}}w_r\overline{w_s}}{\lambda^{\alpha}_r-\lambda^{\alpha}_s}\Big | \leq \frac {\pi}{\alpha\delta} \sum_{r}|w_r|^2.
\end{equation*}
\end{theorem}
  We recall that the Schur product or the Hadamard product of two matrices $X=(x_{i,j})$ and $Y=(y_{i,j})$ is given by  $X \circ Y=(x_{i,j}y_{i,j})$ . A well-known result of Schur asserts that if $X$ is positive definite Hermitian, then $||X \circ Y||_{2,2} \leq \max_i(x_{i,i}) ||Y||_{2,2}$. Schur also showed in \cite{Schur} that the matrix
\begin{equation*}
   X=\Big ( \frac {(\lambda_r-\lambda_s)(\lambda_r\lambda_s)^{\frac {\alpha-1}{2}}}{\lambda^{\alpha}_r-\lambda^{\alpha}_s} \Big )_{r,s}
\end{equation*}
  is positive definite Hermitian for $\alpha \geq 1$, where $X_{r,r}=1/\alpha$. Now, Theorem \ref{thm2} is readily proved following our discussions above applied to $X \circ Y$ with $Y$ being the coefficient matrix of the left-hand side expression of \eqref{2}. In section \ref{sec 4}, we shall give a direct proof of Theorem \ref{thm2}, following the method in \cite{M}. 

\section{A bilinear approach to $l^2$ norms of weighted mean matrices}
\label{sec 2} \setcounter{equation}{0}
   In this section we first recall the following duality principle concerning the norms of linear operators:
\begin{theorem}\cite[Lemma 2]{M}
\label{thm5} Let $C=(c_{n,r})$ be a fixed $N \times R$ matrix. Then the following three assertions concerning the constant $U$ for any ${\bf x} \in l^q, {\bf y} \in l^q$ are equivalent:
\begin{eqnarray*}
 && \Big ( \sum_{r} \Big | \sum_{n} c_{n,r}x_n \Big |^p \Big )^{1/p} \leq  U ||{\bf x}||_p, \\
  && \Big ( \sum_{n} \Big | \sum_{r} c_{n,r}y_n \Big |^q \Big )^{1/q} \leq  U ||{\bf y}||_q, \\
 && \Big | \sum_{n,r} c_{n,r}x_ny_r \Big | \leq  U ||{\bf x}||_p||{\bf y}||_q.
\end{eqnarray*}
\end{theorem}

   For $l^2$ norms of weighted mean matrices, we now point out our general approach. Our goal is to determine, for given $\lambda_k$'s, the best constant $U$ so that the following inequality holds for any ${\bf a} \in l^2$:
\begin{equation}
\label{2.1}
   \sum^{\infty}_{n=1}\Big{|}\sum^{n}_{k=1} \frac {\lambda_ka_k}{\Lambda_n}
   \Big{|}^2 \leq U \sum^{\infty}_{n=1}|a_n|^2.
\end{equation}
   It suffices to consider the cases with the infinite summations above replaced by any finite summations, say from $1$ to $N \geq 1$ here. We may also assume $a_i \geq 0$ for all $i$ here. Now we have
\begin{equation}
\label{2.1'}
   \sum^{N}_{n=1}\Big ( \sum^{n}_{i=1}\frac {\lambda_i}{\Lambda_n}a_i\Big )^2
   =\sum^{N}_{n=1}\Big ( \sum^{n}_{i,j=1}\frac {\lambda_i\lambda_j}{\Lambda^2_n}a_ia_j\Big )
   =\sum^{N}_{i,j=1}\beta_{i,j}a_ia_j, \hspace{0.1in}
   \beta_{i,j}= \sum^{N}_{k = \max{(i,j)}}\frac {\lambda_i\lambda_j}{\Lambda^2_k}.
\end{equation}
  We view the above as a quadratic form and define the associated
  matrix $A$ to be
\begin{equation*}
  A=\Big ( \beta_{i, j} \Big )_{1 \leq i, j \leq N}.
\end{equation*}
   We note that the matrix $A$ is certainly symmetric and positive definite, being equal to $B^{t}B$
   with $B$ a lower-triangular matrix,
\begin{equation*}
   B=\Big ( \eta_{i, j} \Big )_{1 \leq i, j \leq N},  \hspace{0.1in} \eta_{i,j}=\lambda_j/\Lambda_i,  ~~ 1 \leq j \leq
    i;   \hspace{0.1in} \eta_{i,j}=0, j > i.
\end{equation*}
    In order to
   achieve our goal, it suffices to find the maximum eigenvalue of $A$. Instead of working with the matrix $A$ directly, sometimes it is easier to take an equivalent approach. By the duality principle of linear operators (Theorem \ref{thm5} above), the norm of a linear operator acting on the $l^p$ space is equal to the norm of its adjoint acting on the dual space $l^q$. In our case, we may regard inequality  \eqref{2.1} as giving a bound of the norm of the linear operator $B$ acting on the $l^2$ space, which is self-dual. Hence its adjoint $B^t$ (as $B$ is real) also acts on $l^2$ with the same norm as $A$. If we reformulate this in terms of inequalities, what we need to prove are the following Copson (see \cite[Theorems 331, 344]{HLP}) type inequalities:
\begin{equation*}
  \sum^{N}_{n=1}\Big | \sum^{N}_{k=n}\frac {\lambda_n}{\Lambda_k}a_k\Big |^2 \leq U\sum^N_{n=1}|a_n|^2.
\end{equation*}
  Now we have (as $a_i \geq 0$)
\begin{equation}
\label{2.2'}
   \sum^{N}_{n=1}\Big ( \sum^{N}_{k=n}\frac {\lambda_n}{\Lambda_k}a_k\Big )^2
   =\sum^{N}_{n=1}\Big ( \sum^{n}_{i,j=1}\frac {\lambda^2_n}{\Lambda_i\Lambda_j}a_ia_j\Big )
   =\sum^{N}_{i,j=1}\gamma_{i,j}a_ia_j, \hspace{0.1in}
   \gamma_{i,j}= \frac {\sum^{\min (i,j)}_{k=1}\lambda^2_k}{\Lambda_i\Lambda_j}.
\end{equation}

 In summary, in order to establish \eqref{2.1}, we are thus led to the consideration of the following quadratic form (or with $\beta_{i,j}$ replaced by $\gamma_{i,j}$) with ${\bf a}, {\bf b} \in l^2$:
\begin{equation*}
   \sum^{N}_{i,j=1}\beta_{i,j}a_ib_j.
\end{equation*}
   We may assume that $a_i \geq 0, b_j \geq 0$ for any $i, j$ and now we recall the following well-known Schur's test:
\begin{lemma}
\label{lem2.1} Let $p>1$ be fixed and let $A=(\beta_{i,j})_{1 \leq i, j \leq N}$ be a matrix with non-negative entries. If there exist positive numbers $U_1, U_2$ and two positive sequences ${\bf c}=(c_i), 1 \leq i \leq N; {\bf d}=(d_i), 1 \leq i \leq N$, such that
\begin{eqnarray*}
  \sum^{N}_{j=1} \beta_{i,j}c^{1/p}_j  &\leq & U_1d^{1/p}_i,  \hspace{0.1in}  1 \leq i \leq N; \\
 \sum^{N}_{i=1}\beta_{i,j}d^{1/q}_i &\leq & U_2c^{1/q}_j,  \hspace{0.1in}  1 \leq j \leq N.
\end{eqnarray*} 
   Then
\begin{equation*}
  ||A||_{p,p} \leq U^{1/q}_1U^{1/p}_2.
\end{equation*}
\end{lemma}
 
  When we apply the above Lemma to the special case of $p=2, {\bf c}={\bf d}$, we see that (once again by duality) if one can find a positive sequence ${\bf c}$ so that
\begin{equation*}
  \sum^{N}_{j=1} \beta_{i,j}\Big(\frac {c_j}{c_i} \Big)^{1/2}  \leq U_1,  \hspace{0.1in}  1 \leq i \leq N; \hspace{0.1in} 
 \sum^{N}_{i=1}\beta_{i,j}\Big(\frac {c_i}{c_j} \Big )^{1/2} \leq U_2,  \hspace{0.1in}  1 \leq j \leq N,
\end{equation*} 
   for some positive constants $U_1, U_2$, then    
\begin{equation*}
 \sum^{N}_{i,j=1}\beta_{i,j}a_ib_j \leq (U_1U_2)^{1/2}\Big (\sum^{N}_{i=1}a^2_i \Big )^{1/2}\Big ( \sum^{N}_{j=1}b^2_j \Big )^{1/2}.
\end{equation*} 
   From this we see that one can take $U=(U_1U_2)^{1/2}$ in \eqref{2.1}. We point out here that the above approach is the one used in the proof of Theorem 318 in \cite{HLP}. In the special case of $\lambda_k=1$, observe that $\gamma_{i,j}=1/\max (i,j)$ and one takes $c_i=1/i$ to see that on comparing with the corresponding integral (as the integrand is a decreasing function), we have
\begin{equation*}
  \sum^{N}_{j=1} \gamma_{i,j}\Big(\frac {c_j}{c_i} \Big)^{1/2}  \leq \int^N_0\frac {t^{-1/2}}{\max (1, t)}dt=4-\frac {2}{\sqrt{N}}=U_1,
\end{equation*}
  and one can also take $U_2=U_1$ here by symmetry. This implies Hardy's inequality \eqref{eq:1} for $p=2$.
 
  In the next section, we will proceed along the same line above to give a proof of Theorem \ref{thm1}.

\section{Proof of Theorem \ref{thm1}}
\label{sec 3} \setcounter{equation}{0}
   We first look at the case of inequality \eqref{7} for $p=2, 1/2 < \alpha \leq 3/2$. For any real number $\alpha$ satisfying $1/2 < \alpha \leq 3/2$, we denote $M(\alpha)=(m_{i,j})$ for the matrix whose entries are given by
\begin{equation*}
   m_{i,j} =\frac {\alpha^2 \min ( i^{2\alpha-1}, j^{2\alpha-1}) }{(2\alpha-1)i^{\alpha}j^{\alpha}}.
\end{equation*} 
   We apply Lemma \ref{lem2.1} for $c_i=d_i=1/i$ to see that, on comparing with the corresponding integral (as the integrand is a decreasing function)
\begin{eqnarray*}
  \sum^{N}_{j=1} m_{i,j}\Big(\frac {c_j}{c_i} \Big)^{1/2} &= & \frac {\alpha^2}{2\alpha-1} \sum^{N}_{j=1} \frac {1}{i} \frac {(j/i)^{\alpha-3/2}}{\max (1, (j/i)^{2\alpha-1})} \\
& \leq & \frac {\alpha^2}{2\alpha-1} \int^N_0\frac {t^{\alpha-3/2}}{\max (1, t^{2\alpha-1})}dt=\frac {\alpha^2}{(\alpha-1/2)^2}\Big (1- \frac {1}{2\sqrt{N}} \Big )=U_1,
\end{eqnarray*}
   and one can also take $U_2=U_1$ here by symmetry. This now implies the following inequality for $1/2 < \alpha \leq 3/2$ by Lemma \ref{lem2.1}:
\begin{equation}
\label{3.0}
   \sum^{N}_{i,j=1}m_{i,j}a_ia_j= \sum^{N}_{n=1}\Big ( \sum^{N}_{i=n}\frac {\alpha L^{\alpha-1}_{2\alpha-1}(n, n-1)}{i^{\alpha}}a_i\Big )^2 \leq \frac {\alpha^2}{(\alpha-1/2)^2}\sum^{N}_{i=1}a^2_i,
\end{equation}
   for ${\bf a} \in l^2$ with $a_i>0$. Here, the function $L_{2\alpha-1}(n, n-1)$ is defined as in the following lemma with $r=2\alpha-1, a=n, b=n-1$:
\begin{lemma}[{\cite[Lemma 2.1]{alz1.5}}]
\label{lem4}
   Let $a > 0, b > 0$ and $r$ be real numbers with $a \neq b$, and
   let
\begin{eqnarray*}
   L_r (a, b) &=& \genfrac(){1pt}{}{a^r - b^r}{r(a-b)}^{1/(r-1)}  \hspace{0.2in} (r \neq  0,
   1), \\
   L_0(a,b) &=& \frac {a - b}{ \log a-\log
  b}, \\
   L_1(a,b) &=& \frac 1{e}\genfrac(){1pt}{}{a^a}{b^b}^{1/(a-b)}.
\end{eqnarray*}
   The function $r \mapsto L_r(a,b)$ is strictly increasing on ${\mathbb R}$.
\end{lemma}
   
   We now apply the duality principle Theorem \ref{thm5} to \eqref{3.0} to see that it is equivalent to the following inequality:
\begin{equation}
\label{3.3}
   \sum^{N}_{i,j=1}n_{i,j}a_ia_j= \sum^{N}_{n=1}\Big ( \sum^{n}_{i=1}\frac {\alpha L^{\alpha-1}_{2\alpha-1}(i, i-1)}{n^{\alpha}}a_i\Big )^2 \leq \frac {\alpha^2}{(\alpha-1/2)^2}\sum^{N}_{i=1}a^2_i ,
\end{equation}
   for ${\bf a} \in l^2$ with $a_i>0$. Here,
\begin{equation*}
   n_{i,j} =\alpha^2L^{\alpha-1}_{2\alpha-1}(i, i-1)L^{\alpha-1}_{2\alpha-1}(j, j-1) \sum^{N}_{k = \max{(i,j)}}\frac {1}{k^{2\alpha}}.
\end{equation*} 

    Now, to establish \eqref{7}, we just need to show that for every $\alpha>1/2$, $1 \leq i \leq n$,
\begin{equation*}
   \frac {i^{\alpha}-(i-1)^{\alpha}}{n^{\alpha}} \leq \frac {\alpha L^{\alpha-1}_{2\alpha-1}(i, i-1)}{n^{\alpha}},
\end{equation*}
  which is easily seen to follow from Lemma \ref{lem4} with $a=i, b=i-1$ (with $r=\alpha$ and $2\alpha-1$).

  We note that if we denote $N(\alpha)=(n_{i,j})$, then we can also obtain $N(\alpha)$ from $M(\alpha)$ by the following similarity transformation:
\begin{equation*}
  N(\alpha)=FEM(\alpha)E^{-1}F^{-1}, 
\end{equation*}
   where $E=(e_{i,j})_{1 \leq i, j \leq N}, E^{-1}=(e'_{i,j})_{1 \leq i, j \leq N}, E=(f_{i,j})_{1 \leq i, j \leq N}$ are given by
\begin{eqnarray*}
   && e_{i,i}=i^{\alpha}, \hspace{0.1in} e_{i,i-1}=-(i-1)^{\alpha},  \hspace{0.1in} e_{i,j}=0 \hspace{0.1in} \text{otherwise}; \\
&&    e'_{i,j}=\frac {1}{i^{\alpha}}, \hspace{0.1in} j \leq i; \hspace{0.1in} e'_{i,j}=0,  \hspace{0.1in} j > i; \\
   && f_{i,j}=\Big(i^{2\alpha-1}-(i-1)^{2\alpha-1} \Big )^{-1/2}\delta_{i,j}.
\end{eqnarray*} 
 
  We end this section by pointing out that another way of proving \eqref{7} for the case $p=2, 1/2 < \alpha \leq 3/2$ is to show that $\gamma_{i, j} \leq m_{i,j}$ where $\gamma_{i,j}$ is given as in \eqref{2.2'} with $\lambda_k=k^{\alpha}-(k-1)^{\alpha-1}$. In fact, we may assume that $i \leq j$ and observe that
\begin{eqnarray*}
  \gamma_{i,j} &=& \frac {\sum^{\min (i,j)}_{k=1}\lambda^2_k}{\Lambda_i\Lambda_j} =\frac {\sum^{i}_{k=1}\Big ( k^{\alpha}-(k-1)^{\alpha} \Big)^{2}}{i^{\alpha}j^{\alpha}}=\frac {\sum^{i}_{k=1}\Big ( \alpha \int^{k}_{k-1}x^{\alpha-1}dx \Big)^{2}}{i^{\alpha}j^{\alpha}} \\
 & \leq & \frac {\alpha^2 \sum^{i}_{k=1}\int^{k}_{k-1}x^{2(\alpha-1)}dx}{i^{\alpha}j^{\alpha}}=\frac {\alpha^2 i^{2\alpha-1}}{(2\alpha-1)i^{\alpha}j^{\alpha}}=m_{i,j},
\end{eqnarray*}  
  where the inequality above follows from Cauchy's inequality.
 
  Now we consider \eqref{8} for the case $1/p < \alpha \leq 1$, and we show that in this case inequality \eqref{8} follows from the relevant case of inequality \eqref{7}. To show this, we need two lemmas:
\begin{lemma}{\cite[Lemma 2.1]{B5}}
\label{lem2}
    Let ${\bf u}, {\bf v}$ be $n$-tuples
 with non-negative entries with $n \geq 1$ and 
\begin{equation*}
   \sum_{i=1}^k u_i \leq \sum_{i=1}^k v_i, \hspace{0.1in} 1 \leq k \leq n-1; \hspace{0.1in} \sum_{i=1}^n u_i = \sum_{i=1}^n v_i.
\end{equation*}
    then for any decreasing $n$-tuple ${\bf a}$,
\begin{equation*}
   \sum_{i=1}^n u_i a_i \leq \sum_{i=1}^n v_i a_i.
\end{equation*}
\end{lemma}

\begin{lemma}\cite[Lemma 2.4]{CLO}
\label{lem3}
   Let $p>1$ and $A = (a_{n,k})_{n,k \geq 1}$ be a finite lower triangular matrix.
If $a_{n,k} \geq  a_{n,k+1}$ for all $n \geq 1, 1 \leq k \leq n-1$, then
\begin{equation*}
||A||_{p,p} = \sup \Big \{||A{\bf x}||_p: ||{\bf x}||_p = 1 \text{and  ${\bf x}$ decreasing} \Big \}.
\end{equation*}
\end{lemma}
   Now back to the proof of \eqref{8} for the case $1/p < \alpha \leq 1$.   We can replace the infinite sum by a finite sum from $1$ to $N$ with $N \geq 1$ and assume $a_i \geq 0$. It follows that we can apply Lemma \ref{lem3} with the corresponding matrix $A$ being a finite lower triangular matrix with $a_{n,k} \leq  a_{n,k+1}$ for all $n \geq 1, 1 \leq k \leq n-1$ to conclude that it suffices to prove \eqref{8} for the case of ${\bf a}$ being decreasing. It follows from Lemma \ref{lem2} that we will be done if we can show that for any $k \leq n$, 
\begin{equation*}
  \frac {\sum^k_{i=1}i^{\alpha-1}}{\sum^n_{i=1}i^{\alpha-1}} \leq \frac {k^{\alpha}}{n^{\alpha}},
\end{equation*}
   since one knows already that \eqref{7} holds. By induction, it suffices to show the above inequalities for the case $k=n-1$. In this case we can recast the above inequality as
\begin{equation*}
  P_{n-1}(\alpha-1) \geq \frac {n-1}{n},
\end{equation*}
 where for any integer $n \geq 1$ and any real number $r>0$, we define
\begin{equation*}
  P_{n}(r)= \left(\frac {1}{n} \sum_{i=1}^{n}i^r\bigg/
\frac {1}{n+1}\sum_{i=1}^{n+1}i^r\right)^{1/r}.
\end{equation*}
   Our proof now follows from a combination of two results of Alzer \cite{alz} and \cite[Theorem 2.3]{alz1}, which asserts that for any real $r$,
\begin{equation*}
   P_n(r) > \frac {n}{n+1}= \lim_{r \rightarrow +\infty}P_n(r).
\end{equation*}

\section{Proof of Theorem \ref{thm2}}
\label{sec 4} \setcounter{equation}{0}
    We now give the proof of Theorem \ref{thm2} and we reserve the letter $i$ for the complex number $\sqrt{-1}$ in this section.  From the duality principle of linear operators Theorem \ref{thm5}, it suffices to show that
\begin{equation*}
   \sum_r \Big |
\sum_{\substack{s \\s \neq r}} \frac {(\lambda_r\lambda_s)^{\frac {\alpha-1}{2}}\overline{w_s}}{\lambda^{\alpha}_r-\lambda^{\alpha}_s}\Big |^2 \leq \frac {\pi^2}{\alpha^2\delta^2} \sum_{r}|w_r|^2.
\end{equation*}
   We multiply out the square on the left and take the sum over $r$ inside to see that the left-hand side is
\begin{equation*}
   \sum_{s,t} \overline{w_s}w_t
\sum_{\substack{r \\r \neq s,t}} \frac {(\lambda_r\lambda_s)^{\frac {\alpha-1}{2}}}{\lambda^{\alpha}_r-\lambda^{\alpha}_s} \cdot \frac {(\lambda_r\lambda_t)^{\frac {\alpha-1}{2}}}{\lambda^{\alpha}_r-\lambda^{\alpha}_t}. 
\end{equation*}
   Writing the diagonal terms separately, we see that this is
\begin{equation*}
  \sum_{s} |w_s|^2
\sum_{\substack{r \\r \neq s}} \frac {(\lambda_r\lambda_s)^{\alpha-1}}{(\lambda^{\alpha}_r-\lambda^{\alpha}_s)^2} + \sum_{\substack{s,t \\s \neq t}} \overline{w_s}w_t  \frac {(\lambda_s\lambda_t)^{\frac {\alpha-1}{2}}}{\lambda^{\alpha}_s-\lambda^{\alpha}_t}
\sum_{\substack{r \\r \neq s,t}} \Big ( \frac {\lambda_r^{\alpha-1}}{\lambda^{\alpha}_r-\lambda^{\alpha}_s} - \frac {\lambda_r^{\alpha-1}}{\lambda^{\alpha}_r-\lambda^{\alpha}_t} \Big ). 
\end{equation*}
  For the second term above, we write it as the difference of two terms in which the inner summands are $\lambda_r^{\alpha-1}(\lambda^{\alpha}_r-\lambda^{\alpha}_s)^{-1}$ and $\lambda_r^{\alpha-1}(\lambda^{\alpha}_r-\lambda^{\alpha}_t)^{-1}$, respectively. In the first of these we introduce the new term for $r=t$, and similarly for the second. Thus the second term above is
\begin{eqnarray*}
  &=& \sum_{\substack{s,t \\s \neq t}} \overline{w_s}w_t  \frac {(\lambda_s\lambda_t)^{\frac {\alpha-1}{2}}}{\lambda^{\alpha}_s-\lambda^{\alpha}_t}
\sum_{\substack{r \\r \neq s}} \frac {\lambda_r^{\alpha-1}}{\lambda^{\alpha}_r-\lambda^{\alpha}_s}-\sum_{\substack{s,t \\s \neq t}} \overline{w_s}w_t  \frac {(\lambda_s\lambda_t)^{\frac {\alpha-1}{2}}}{\lambda^{\alpha}_s-\lambda^{\alpha}_t}
\sum_{\substack{r \\r \neq t}} \frac {\lambda_r^{\alpha-1}}{\lambda^{\alpha}_r-\lambda^{\alpha}_t}  \\
&& +  \sum_{\substack{s,t \\s \neq t}} \overline{w_s}w_t  \frac {(\lambda_s\lambda_t)^{\frac {\alpha-1}{2}}\lambda^{\alpha-1}_t}{(\lambda^{\alpha}_s-\lambda^{\alpha}_t)^2}+  \sum_{\substack{s,t \\s \neq t}} \overline{w_s}w_t  \frac {(\lambda_s\lambda_t)^{\frac {\alpha-1}{2}}\lambda^{\alpha-1}_s}{(\lambda^{\alpha}_s-\lambda^{\alpha}_t)^2}.
\end{eqnarray*}
  We denote $\sum_1$ and $\sum_2$ for the first two terms above, respectively. Now we may assume that the $w_r$'s are extreme and as the coefficient matrix is skew-Hermitian, the extremal $w_r$'s are the coordinates of an eigenvector. Hence, there is a real number $\mu$ such that
\begin{equation*}
  \sum_{\substack{r \\r \neq s}} \frac {(\lambda_r\lambda_s)^{\frac {\alpha-1}{2}}w_r}{\lambda^{\alpha}_r-\lambda^{\alpha}_s}=i\mu w_s, \hspace{0.1in} 1 \leq s \leq R.
\end{equation*}
  Taking the sum over $t$ inside in $\sum_1$, and using the above, we find that
\begin{equation*}
 \sum_{\substack{s,t \\s \neq t}} \overline{w_s}w_t  \frac {(\lambda_s\lambda_t)^{\frac {\alpha-1}{2}}}{\lambda^{\alpha}_s-\lambda^{\alpha}_t}
\sum_{\substack{r \\r \neq s}} \frac {\lambda_r^{\alpha-1}}{\lambda^{\alpha}_r-\lambda^{\alpha}_s}=-i\mu\sum_{s} |w_s|^2 
\sum_{\substack{r \\r \neq s}} \frac {\lambda_r^{\alpha-1}}{\lambda^{\alpha}_r-\lambda^{\alpha}_s}.
\end{equation*}
  Making the same simplification in $\sum_2$ we find that $\sum_1=\sum_2$ for those extreme $w_r$'s. Thus we have
\begin{eqnarray*}
 &&  \sum_r \Big |
\sum_{\substack{s \\s \neq r}} \frac {(\lambda_r\lambda_s)^{\frac {\alpha-1}{2}}\overline{w_s}}{\lambda^{\alpha}_r-\lambda^{\alpha}_s}\Big |^2 \\
 & \leq & \sum_{s} |w_s|^2
\sum_{\substack{r \\r \neq s}} \frac {(\lambda_r\lambda_s)^{\alpha-1}}{(\lambda^{\alpha}_r-\lambda^{\alpha}_s)^2}+  \sum_{\substack{s,t \\s \neq t}} \overline{w_s}w_t  \frac {(\lambda_s\lambda_t)^{\frac {\alpha-1}{2}}(\lambda^{\alpha-1}_s+\lambda^{\alpha-1}_t)}{(\lambda^{\alpha}_s-\lambda^{\alpha}_t)^2} \\
& \leq & \sum_{s} |w_s|^2
\sum_{\substack{r \\r \neq s}} \frac {(\lambda_r\lambda_s)^{\alpha-1}}{(\lambda^{\alpha}_r-\lambda^{\alpha}_s)^2}+  \sum_{\substack{s,t \\s \neq t}} \Big ( \frac {|w_s|^2}{2}+\frac {|w_t|^2}{2} \Big )  \frac {(\lambda_s\lambda_t)^{\frac {\alpha-1}{2}}(\lambda^{\alpha-1}_s+\lambda^{\alpha-1}_t)}{(\lambda^{\alpha}_s-\lambda^{\alpha}_t)^2} \\
&=& \sum_{s} |w_s|^2
\sum_{\substack{r \\r \neq s}} \frac {(\lambda_r\lambda_s)^{\alpha-1}+(\lambda_r\lambda_s)^{\frac {\alpha-1}{2}}(\lambda^{\alpha-1}_r+\lambda^{\alpha-1}_s)}{(\lambda^{\alpha}_r-\lambda^{\alpha}_s)^2},
\end{eqnarray*}
  where we have used $|\overline{w_s}w_t| \leq \frac {|w_s|^2}{2}+\frac {|w_t|^2}{2}$ in the second inequality above. To establish Theorem \ref{thm2}, it now suffices to show that for $\alpha \geq 1$,
\begin{equation}
\label{1.3}
   \frac {(\lambda_r\lambda_s)^{\alpha-1}+(\lambda_r\lambda_s)^{\frac {\alpha-1}{2}}(\lambda^{\alpha-1}_r+\lambda^{\alpha-1}_s)}{(\lambda^{\alpha}_r-\lambda^{\alpha}_s)^2} \leq \frac {3}{\alpha^2(\lambda_r-\lambda_s)^2},
\end{equation}
  as the rest will follow from the treatment in the proof of Theorem 2 in \cite{M}.

   To establish \eqref{1.3} for $\alpha \geq 1$, we note first that by the arithmetic-geometric mean inequality, 
\begin{equation*}
  (\lambda_r\lambda_s)^{\frac {\alpha-1}{2}} \leq \frac {\lambda^{\alpha-1}_r+\lambda^{\alpha-1}_s}{2}.
\end{equation*}
   Hence it suffices to show that
\begin{equation*}
   \frac {(\lambda_r\lambda_s)^{\frac {\alpha-1}{2}}(\lambda^{\alpha-1}_r+\lambda^{\alpha-1}_s)}{(\lambda^{\alpha}_r-\lambda^{\alpha}_s)^2} \leq \frac {2}{\alpha^2(\lambda_r-\lambda_s)^2}.
\end{equation*}
  Without loss of generality, we may assume $\lambda_r \geq \lambda_s$ and on letting $t=\lambda_r / \lambda_s$, we can recast the above inequality as
\begin{equation*}
   \frac {t^{\frac {\alpha-1}{2}}(t^{\alpha-1}+1)}{2} \leq \Big ( \frac {t^{\alpha}-1}{\alpha(t-1)} \Big )^2.
\end{equation*}
   We now set $y=t^{\alpha-1}$ to rewrite the above inequality as
\begin{equation}
\label{1.4}
   \frac {y^{\frac {1}{2}}(y+1)}{2} \leq \Big ( \frac {y^{\alpha/(\alpha-1)}-1}{\alpha(y^{1/(\alpha-1)}-1)} \Big )^2.
\end{equation}
   We let $\beta=1/(\alpha-1)$ to write the right-hand side expression above as
\begin{equation*}
    \frac {f^2(1+\beta)}{f^2(\beta)} , \hspace{0.1in} f(\beta)=\frac {y^{\beta}-1}{\beta}.
\end{equation*} 
   We want to show the above function is increasing with $\beta$ for $\beta>0$, which one checks is equivalent to
\begin{equation*}
    \frac {f'(1+\beta)}{f(1+\beta)} \geq \frac {f'(\beta)}{f(\beta)}.
\end{equation*} 
   The above inequality follows if one can show that $f'(\beta)/f(\beta)$ is an increasing function of $\beta$. Direct calculations show that
\begin{equation*}
    \Big ( \frac {f'(\beta)}{f(\beta)} \Big )'=-\frac {\ln^2 y \cdot y^{\beta}}{(y^{\beta}-1)^2}+\frac {1}{\beta^2}.
\end{equation*} 
   To show the right-hand side expression above is $\geq 0$, it is equivalent to showing that (note that $y \geq 1$ here)
\begin{equation*}
   \frac {y^{\beta}-1}{\beta}=\frac {\int^{\beta}_0 \ln y \cdot y^s ds}{\beta-0} \geq \ln y \cdot y^{\beta/2}.
\end{equation*} 
   The above inequality now follows from the well-known Hadamard's inequality (with $h(x)=\ln y \cdot y^{x}, a=0, b=\beta$ here), which asserts that for a continuous convex function $h(x)$ on $[a, b]$,
\begin{equation}
\label{3.6}
   h(\frac {a+b}2) \leq \frac {1}{b-a}\int^b_a h(x)dx \leq \frac
   {h(a)+h(b)}{2}.
\end{equation}
   It follows now that 
\begin{equation*}
    \Big ( \frac {f(1+\beta)}{f(\beta)} \Big )^2 \geq \lim_{\beta \rightarrow 0^+}\Big ( \frac {f(1+\beta)}{f(\beta)} \Big )^2=\Big ( \frac {y-1}{\ln y} \Big )^2.
\end{equation*} 
   Thus in order to prove \eqref{1.4}, it suffices to show that
\begin{equation*}
     \frac {y^{\frac {1}{2}}(y+1)}{2} \leq \Big ( \frac {y-1}{\ln y} \Big )^2.
\end{equation*} 
   We rewrite the above as
\begin{equation*}
     \Big ( \frac {\ln y-\ln 1}{y-1} \Big )^2=\Big ( \frac {\int^y_1 \ln s ds}{y-1} \Big )^2 \leq \frac {2}{y^{\frac {1}{2}}(y+1)} .
\end{equation*} 
   Apply Hadamard's inequality \eqref{3.6} again here with $h(x)=-\ln x , a=0, b=y$, we see that
\begin{equation*}
    \Big ( \frac {\int^y_1 \ln s ds}{y-1} \Big )^2 \leq \Big ( \frac {2}{y+1} \Big )^2 \leq \frac {2}{y^{\frac {1}{2}}(y+1)},
\end{equation*} 
   where the last inequality follows easily from the arithmetic-geometric mean inequality and this completes the proof of Theorem \ref{thm2}.

\section{Some Remarks}
\label{sec 5} \setcounter{equation}{0}
  We point out here that the matrix $(\gamma_{i,j})$ given in \eqref{2.2'} is related to the matrix $(\beta_{i,j})$ as given in \eqref{2.1'}  by the following similarity transformation:
\begin{equation*}
  \Big ( \gamma_{i, j} \Big )_{1 \leq i, j \leq N}=GHAH^{-1}G^{-1}, 
\end{equation*}
   where $A=\Big ( \beta_{i, j} \Big )_{1 \leq i, j \leq N}, H=(h_{i,j})_{1 \leq i, j \leq N}, H^{-1}=(h'_{i,j})_{1 \leq i, j \leq N}, G=(g_{i,j})_{1 \leq i, j \leq N}$ with
\begin{eqnarray*}
   && h_{i,j}=\lambda_j, \hspace{0.1in} j \leq i; \hspace{0.1in} h_{i,j}=0,  \hspace{0.1in} j > i; \\
   && h'_{i,i}=1/\lambda_i, \hspace{0.1in} h'_{i,i-1}=-1/\lambda_i,  \hspace{0.1in} h_{i,j}=0 \hspace{0.1in} \text{otherwise}; \\
   && g_{i,j}=\frac {\Lambda_1}{\Lambda_i}\delta_{i,j}.
\end{eqnarray*} 
    We note here in the special case of $\lambda_k=1$, which corresponds to Hardy's inequality \eqref{eq:1}, the above relation is used in \cite{W1}, in which case $\gamma_{i,j}=1/\max (i,j)$. As $(\gamma_{i,j})$ and $A$ have the same set of eigenvalues, we are again led to our approach above.
 
  We note that if one can show that for any $i,j$,
\begin{equation}
\label{2.3}
  \gamma_{i,j} \geq 1/\max(i, j),
\end{equation}
   then inequality \eqref{2.1} holds with $U \geq 4$, as the right-hand side expression above corresponds to the value of $\gamma_{i,j}$ when $\lambda_k=1$. We may now assume that $i \leq j$ and note that \eqref{2.3} holds when $i=j$ by Cauchy's inequality. For the general case, note that
\begin{equation*}
  \gamma_{i,j} = \frac {\sum^{\min (i,j)}_{k=1}\lambda^2_k}{\Lambda_i\Lambda_j}= \frac {\sum^{\min (i,j)}_{k=1}\lambda^2_k}{\Lambda^2_i}\frac {\Lambda_i}{\Lambda_j} \geq \frac {\Lambda_i}{i\Lambda_j}.
\end{equation*}
  It is easy to see that when the $\lambda_k$'s are decreasing, then the right-hand side expression above is no less than the right-hand side expression of \eqref{2.3} and hence we have just shown that when the $\lambda_k$'s are decreasing, then inequality \eqref{2.1} holds with $U \geq 4$.
This now gives another proof of a result of Bennett \cite[Theorem 1.14 (a)]{B5} for the case of weighted mean matrices when $p=2$. We point out here that Bennett's result \cite[Theorem 1.14 ]{B5} is more general and covers the case for any summability matrix. The non-trivial part of his result \cite[Theorem 1.14(b)]{B5}, which asserts that for any summability matrix $A$, if the rows of $A$ are increasing, then 
\begin{equation*}
  ||A||_{p,p} \leq q,
\end{equation*} 
   is a consequence of \cite[Theorem 332]{HLP} when $A$ is a weighted mean matrix. In fact, by a change of variables $a_n \rightarrow \lambda^{-1/p}_na_n$, we see that it yields
\begin{equation*}
  \sum^{\infty}_{n=1}\Big{|}\frac {1}{\Lambda_n}
\sum^n_{k=1}\lambda^{1/q}_k\lambda^{1/p}_na_k\Big{|}^p \leq \Big (\frac
{p}{p-1} \Big )^p\sum^\infty_{k=1}|a_k|^p,
\end{equation*}
  from which one deduces Bennett's result easily.

   
   We end this paper by pointing out that in view of \eqref{3.3} and Lemma \ref{lem4}, the following inequality holds for any ${\bf a} \in l^2$ for $1/2 < \alpha  \leq 3/2, \max (2\alpha-1, \alpha) \geq s \geq \min (2\alpha -1, \alpha)$:
\begin{equation*}
   \sum^{N}_{n=1}\Big | \sum^{n}_{i=1}\frac {\alpha L^{\alpha-1}_{s}(i, i-1)}{n^{\alpha}}a_i\Big|^2 \leq \frac {\alpha^2}{(\alpha-1/2)^2}\sum^{N}_{i=1}|a_i|^2.
\end{equation*}


\end{document}